\documentclass[12pt]{article}
\usepackage{amsmath}
\usepackage{amssymb}
\usepackage{amsmath,amssymb,amsbsy,amsfonts,amsthm,latexsym,
amsopn,amstext,amsxtra,euscript,amscd}

\newtheorem{lem}{Lemma}
\newtheorem{lemma}[lem]{Lemma}
\newtheorem{thm}{Theorem}
\newtheorem{theorem}[thm]{Theorem}

\begin{document}

\baselineskip 15pt

\title{An explicit sum-product estimate in $\mathbb{F}_p$}

\author{{M. Z. Garaev}
\\
\normalsize{Instituto de Matem{\'a}ticas,  UNAM}
\\
\normalsize{Campus Morelia, Apartado Postal 61-3 (Xangari)}
\\
\normalsize{C.P. 58089, Morelia, Michoac{\'a}n, M{\'e}xico} \\
\normalsize{\tt garaev@matmor.unam.mx}\\
}

\date{\empty}

\pagenumbering{arabic}

\maketitle

\begin{abstract}
Let $\mathbb{F}_p$ be the field of residue classes modulo a prime
number $p$ and let $A$ be a non-empty subset of $\mathbb{F}_p.$  In
this paper we give an explicit version of the sum-product estimate
of Bourgain, Katz, Tao and Bourgain, Glibichuk, Konyagin on the size
of $\max\{|A+A|, |AA|\}.$ In particular, our result implies that if
$1<|A|\le p^{7/13}(\log p)^{-4/13},$ then
$$
\max\{|A+A|, |AA|\}\gg \frac{|A|^{15/14}}{(\log|A|)^{2/7} } .
$$
\end{abstract}



\section{Introduction}

Let $\mathbb{F}_p$ be the field of residue classes modulo a prime
number $p$ and let  $A$ be a non-empty subset of $\mathbb{F}_p.$
Consider the sum set
$$
A+A=\{a+b: \, a\in A, \, b\in A\}
$$
and the product set
$$
AA =\{ab: \, a\in A, \, b\in A\}.
$$
Bourgain, Katz, Tao~\cite{BKT} and Bourgain, Glibichuk,
Konyagin~\cite{BGK} have shown that if $|A|<p^{1-\delta},$ where
$\delta>0,$ then one has the sum-product estimate
\begin{equation}
\label{eqn:sumprodclass} \max\{|A+A|,|AA|\}\gg |A|^{1+\varepsilon}
\end{equation}
for some $\varepsilon=\varepsilon(\delta)>0.$ This result has found
a number of spectacular applications in combinatorial problems and
exponential sum estimates, see~\cite{B1}--\cite{BKT}.

Bound~\eqref{eqn:sumprodclass} does not yield an explicit
relationship between $\varepsilon$ and $\delta.$ Hart, Iosevich and
Solymosi~\cite{HIS}, by using Kloosterman sums, could obtain a
concrete value of $\varepsilon$ in certain ranges of $|A|.$ More
precisely, they proved that
\begin{equation}
\label{eqn:HIS} \max\{|A+A|,\, |AA|\}\gg \left\{ \begin{array}{ll}
|A|^{3/2}p^{-1/4}, & \quad \mbox{if}\quad
p^{1/2}<|A|<p^{7/10},\\
|A|^{2/3}p^{1/3}, & \quad \mbox{if}\quad p^{7/10}<|A|\le p.
\end{array} \right.
\end{equation}
The aim of the present paper is to obtain an explicit sum-product
estimate for any range of $|A|.$
\begin{theorem}
\label{thm:sumprodest} Let $|A|>1.$ Then
$$
\max\{|A+A|, |AA|\}\gg \min\left\{\frac{|A|^{15/14}\max\Bigl\{1,\,
|A|^{1/7}p^{-1/14}\Bigr\}}{(\log |A|)^{2/7}},\,
\frac{|A|^{11/12}p^{1/12}}{ (\log|A|)^{1/3}}\right\}.
$$
In particular, if $1<|A|< p^{7/13}(\log p)^{-4/13},$ then we have
$$
\max\{|A+A|, |AA|\}\gg \frac{|A|^{15/14}}{(\log|A|)^{2/7} } .
$$
\end{theorem}

The proof of Theorem~\ref{thm:sumprodest} uses results and tools
from arithmetical combinatorics. When $|A|$ is larger than
$p^{5/9}$, combinatorial arguments and trigonometric sums can be
used together to get a better estimate.
\begin{theorem}
\label{thm:sumprodest1} Let $|A|>1.$ Then
$$
\max\{|A+A|, |AA|\}\gg \min\Bigl\{|A|^{5/3}p^{-1/3}(\log
|A|)^{-1/3}, |A|^{2/3}p^{1/3} (\log|A|)^{-1/3}\Bigr\}.
$$
In particular, if $1<|A|<p^{2/3},$ then we have
$$
\max\{|A+A|, |AA|\}\gg |A|^{5/3}p^{-1/3}(\log |A|)^{-1/3}.
$$
\end{theorem}
When $p^{7/10}(\log p)^{-1/3}<|A|<p,$ the inequality~\eqref{eqn:HIS}
is preferable. We remark that~\eqref{eqn:HIS} can also be proved
using the method described in the present paper.

In the corresponding problem for integers (i.e., if the field
$\mathbb{F}_p$ is replaced by the set of integers) the conjecture of
Erd\"{o}s and Szemer\'edi~\cite{ErSz} is that $\max\{|A+A|,|AA|\}\gg
c(\varepsilon)|A|^{2-\varepsilon}$ for any $\varepsilon>0.$ At
present the best known bound in the integer problem is
$\max\{|A+A|,|AA|\}\gg |A|^{14/11}(\log|A|)^{-3/11}$ due to
Solymosi~\cite{Sol}.

We do not know what the optimal lower bound for $\max\{|A+A|,\,
|AA|\}$ in terms of $|A|$ and $p$ should be. It is known that the
analogy of the Erd\"{o}s and Szemer\'edi conjecture in the form
$\max\{|A+A|,|AA|\}\ge c(\varepsilon)\min\{|A|^{2-\varepsilon},
p^{1-\varepsilon}\}$ for all subsets $A$ of $\mathbb{F}_p$ does not
hold.

In what follows, all the sets under consideration are assumed to be
non-empty. For a set $X\subset \mathbb{F}_p$ and for an element
$a\in \mathbb{F}_p$ we use the notation
$$
a*X=\{ax:\, x\in X\}.
$$

\section{Lemmas}

The following lemma follows from the work of Glibichuk and
Konyagin~\cite{GK}.

\begin{lemma}
\label{lem:GKon} Let $A_1\subset \mathbb{F}_p$ with $1<|A_1|<
p^{1/2}.$ Then there exist elements $a_1, a_2, b_1, b_2\in A_1$ such
that $a_1\not=a_2$ and
$$
|(a_1-a_2)*A_1+(a_1-a_2)*A_1+(b_1-b_2)*A_1|\ge 0.5|A_1|^2.
$$
\end{lemma}

When $|A_1|>p^{1/2},$ we will use the following statement
(see~\cite{BGK} or~\cite{Glib}).
\begin{lemma}
\label{lem:bigAgeneral} Let $A_1\subset \mathbb{F}_p$ with $|A_1|>
p^{1/2}.$ Then there exist  $a_1, a_2, b_1, b_2 \in A_1$ such that
$$
|(a_1-a_2)*A_1+(b_1-b_2)*A_1|\ge 0.5p.
$$
\end{lemma}

The following lemmas are due to Ruzsa
(see~\cite{Nat},~\cite{R1},~\cite{R2},~\cite{TV}). They hold for
subsets of any abelian group, but here we state them only for
subsets of $\mathbb{F}_p.$ We will repeatedly use these lemmas to
prove our result. Lemma~\ref{lem:RuzsaTr} is called Ruzsa's triangle
inequality.

\begin{lemma}
\label{lem:RuzsaTr} For any subsets $X,Y,Z$ of $\mathbb{F}_p$ we
have
$$
|X-Z|\le\frac{|X-Y||Y-Z|}{|Y|}.
$$
\end{lemma}

\begin{lemma}
\label{lem:RuzsaPlun} Let $X,B_1,\ldots, B_k$ be any subsets of
$\mathbb{F}_p$ with
$$
|X|=n,\quad |X+B_i|=\alpha_i n, \quad (i=1, \ldots, k).
$$
Then there is an $X_1\subset X$ such that
$$
|X_1+B_1+\ldots+B_k|\le \alpha_1\ldots \alpha_k |X_1|
$$
\end{lemma}

An important corollary of Lemma~\ref{lem:RuzsaPlun} is the
inequality
$$
|B_1+\ldots+B_k|\le \frac{|X+B_1|\ldots |X+B_{k}|}{|X|^{k-1}}.
$$
Below, when we refer to Lemma~\ref{lem:RuzsaPlun}, we will always
understand this inequality. The case $k=2$ illustrates another
version of Ruzsa's triangle inequality.

\section{Proof of Theorem~\ref{thm:sumprodest}}

We use an idea of the proof of Katz-Tao lemma presented
in~\cite[Section 2.8]{TV}. That proof, as it was mentioned
in~\cite{TV}, used Bourgain's idea from~\cite{B2}.

We may assume that $|A|^2\ge 100|AA|$ and that $0\not\in A.$ Let $J$
denote the number of solutions of the equation
$$
ax=by,\quad a,b,x,y\in A.
$$
From the well-known relationship between the cardinality of a set
and the number of solutions of the corresponding equation, we have
that
$$
J\ge \frac{|A|^4}{|AA|}.
$$
Since
$$
J=\sum_{a\in A}\sum_{b\in A}|a*A\cap b*A|,
$$
there exists a fixed element $b=b_0\in A$ such that
$$
\sum_{a\in A}|a*A\cap b_0*A|\ge \frac{|A|^3}{|AA|}.
$$
For a given positive integer $j\le \log|A|/\log 2 +1,$ let $D_j$ be
the set of all $a\in A$ for which
$$
2^{j-1}\le |a*A\cap b_0*A|< 2^j.
$$
Then,
$$
\sum_j\sum_{a\in D_j} 2^j\ge \frac{|A|^3}{|AA|}.
$$
Let the quantity $\sum\limits_{a\in D_j} 2^j$ takes its biggest
value (or one of its biggest values if there are several such ones)
when  $j=j_0.$ Denote
$$
N=2^{j_0-1},\qquad A_1=D_{j_0}\subset A.
$$ 
Then, for any $a\in A_1,$ we have
\begin{equation}
\label{eqn:aAcapbAge} 1\le N\le |a*A\cap b_0*A|< 2N, \quad
N|A_1|\gg\frac{|A|^3}{|AA|\log|A|}.
\end{equation}
In particular, since $N\le |A|$ and $|A_1|\le |A|,$ we get
\begin{equation}
\label{eqn:boundA1} |N|\gg\frac{|A|^2}{|AA|\log|A|},\qquad
|A_1|\gg\frac{|A|^2}{|AA|\log|A|}.
\end{equation}
Now let $a$ be an arbitrary element of $A_1.$ From
Lemma~\ref{lem:RuzsaTr} and the inequality~\eqref{eqn:aAcapbAge} we
have
\begin{eqnarray*}
&& |a*A-b_0*A|\le \\
&& \frac{|a*A+(a*A\cap b_0*A)||-(a*A\cap b_0*A)-b_0*A|}{{|a*A\cap
b_0*A|}}\le \\
&& \frac{|a*A+a*A||b_0*A+b_0*A|}{N}= \frac{|A+A|^2}{N}.
\end{eqnarray*}
Furthermore, using Lemma~\ref{lem:RuzsaPlun} with
$$
k=2,\,\, X=a*A\cap b_0*A,\,\, B_1=a*A, \,\,B_2=b_0*A,
$$
we obtain
\begin{eqnarray*}
&& |a*A+b_0*A|\le \frac{|a*A+(a*A\cap b_0*A)||b_0*A+(a*A\cap
b_0*A)|}{|a*A\cap b_0*A|}\le\\
&&\frac{|a*A+a*A||b_0*A+b_0*A|}{N}=\frac{|A+A|^2}{N}.
\end{eqnarray*}
Thus, the bound
\begin{equation}
\label{eqn:aA pm b0A} |a*A\pm b_0*A|\ll \frac{|A+A|^2}{N}
\end{equation}
holds for any $a\in A_1$ and for any choice of the sign ``$\pm$".

There are two cases to consider.

\bigskip

{\bf Case 1:} $|A_1|< p^{1/2}.$

\bigskip

First of all, in this case besides of~\eqref{eqn:boundA1}, we also
have
$$
Np^{1/2}>N|A_1|\ge\frac{|A|^3}{|AA|\log|A|}.
$$
Thus, together with~\eqref{eqn:boundA1}, we have
\begin{equation}
\label{eqn:Case1NA1} N>\max\Bigl\{\frac{|A|^2}{|AA|\log|A|}, \,
\frac{|A|^3}{p^{1/2}|AA|\log|A|}\Bigr\}.
\end{equation}
 Since $A_1\subset A,$ according to
Lemma~\ref{lem:GKon} there exist $a_1,a_2,b_1,b_2\in A_1$ such that
$a_1\not=a_2$ and
$$
0.5|A_1|^2 \le |(a_1-a_2)*A+(a_1-a_2)*A+(b_1-b_2)*A|.
$$
We apply Lemma~\ref{lem:RuzsaPlun} with $k=3$ and
$$
X=B_1=B_2=(a_1-a_2)*A, \quad B_3=(b_1-b_2)*A.
$$
Then we get
\begin{equation}
\label{eqn:3A-3A to 2A-2A} 0.5|A_1|^2 \le
\frac{|A+A|^{2}\Bigl|(a_1-a_2)*A+(b_1-b_2)*A\Bigr|}{|A|^2}.
\end{equation}
Next, we apply Lemma~\ref{lem:RuzsaPlun} with $k=4$ and
$$
X=b_0*A, \,\, B_1=a_1*A, \,\, B_2=-a_2*A, \,\, B_3=b_1*A, \,\,
B_4=-b_2*A.
$$
Then,
\begin{eqnarray}
\label{eqn:2A-2A to b0A-aA}
\begin{split}
& |(a_1-a_2)*A+(b_1-b_2)*A|\le\\
& \quad
\frac{|b_0*A+a_1*A||b_0*A-a_2*A||b_0*A+b_1*A||b_0*A-b_2*A|}{|A|^3}.
\end{split}
\end{eqnarray}
Applying the inequality~\eqref{eqn:aA pm b0A} to the right hand side
of~\eqref{eqn:2A-2A to b0A-aA} and incorporating the resulting
estimate to~\eqref{eqn:3A-3A to 2A-2A}, we get
$$
|A+A|^{10}\gg |A_1|^2|A|^5N^4.
$$
Taking into account~\eqref{eqn:aAcapbAge} to substitute $N|A_1|$ and
then~\eqref{eqn:Case1NA1} to substitute $N,$ we conclude that
$$
|A+A|^{10}|AA|^4(\log|A|)^4\gg |A|^{15}\max\Bigl\{1,
\frac{|A|^2}{p}\Bigr\}.
$$
This proves Theorem~\ref{thm:sumprodest} in Case 1.

\bigskip

{\bf Case 2}:\quad $|A_1|>p^{1/2}.$

\bigskip

In this case, according to Lemma~\ref{lem:bigAgeneral}, there exist
elements $a_1,a_2,b_1,b_2\in A_1$ such that
$$
0.5p \le |(a_1-a_2)*A+(b_1-b_2)*A|.
$$
Applying Lemma~\ref{lem:RuzsaPlun} with $k=4$ and
$$
X=b_0*A,\,\, B_1=a_1*A,\,\, B_2=-a_2*A, \,\, B_3=b_1*A,
\,\,B_4=-b_2*A,
$$
we get
$$
0.5p\le
\frac{|b_0*A+a_1*A||b_0*A-a_2*A||b_0*A+b_1*A||b_0*A-b_2*A|}{|A|^3}.
$$
Taking into account the inequality~\eqref{eqn:aA pm b0A}, we obtain
that
$$
p\ll \frac{|A+A|^{8}}{|A|^{3}N^4}.
$$
Using~\eqref{eqn:boundA1}, we get
$$
|A+A|^8|AA|^4(\log|A|)^4\gg p|A|^{11}.
$$
Theorem~\ref{thm:sumprodest} is proved.

\section{Proof of Theorem~\ref{thm:sumprodest1}}

We remark that the proof of Lemma~\ref{lem:bigAgeneral}, as well as
the proof of Lemma~\ref{lem:GKon}, uses the fact that for any sets
$X, Y, G\subset \mathbb{F}_p$ there exists $\xi\in G$ such that
$$
|X+\xi*Y|\ge \frac{|X||Y||G|}{|X||Y|+|G|}.
$$
This estimate is nontrivial when $|G|$ is larger than $|X|$ and
$|Y|.$ In order to prove Theorem~\ref{thm:sumprodest1} we would like
to have a similar estimate which in certain cases would be
nontrivial when the cardinality of $G$ is smaller than those of $X$
and $Y.$  The following lemma provides with such an estimate.
\begin{lemma}
\label{lem:XYG} Let $X, Y, G\subset \mathbb{F}_p.$ Then there exists
$\xi\in G$ such that
$$
|X+\xi*Y|\ge \frac{p|X||Y||G|}{|X||Y||G|+p^2}.
$$
\end{lemma}

It is easy to see that the
bounds~\eqref{eqn:aAcapbAge},~\eqref{eqn:boundA1},~\eqref{eqn:aA pm
b0A} and Lemma~\ref{lem:XYG} imply Theorem~\ref{thm:sumprodest1}.
Indeed, apply Lemma~\ref{lem:XYG} with
$$
X=-b_0*A,\quad Y=A,\quad G=A_1.
$$
Then for some $a\in A_1,$ in view of~\eqref{eqn:aA pm b0A}, we have
$$
\frac{|A+A|^2}{N}\gg |a*A-b_0*A|\gg \min\Bigl\{p, \,
\frac{|A|^2|A_1|}{p}\Bigr\}.
$$
Thus,
$$
|A+A|^2\gg \min\Bigl\{pN, \, \frac{|A|^2N|A_1|}{p}\Bigr\}\gg
\min\Bigl\{p|A|^2|AA|^{-1}, \, \frac{|A|^2N|A_1|}{p}\Bigr\}.
$$
Recalling~\eqref{eqn:aAcapbAge},~\eqref{eqn:boundA1}, we get that
either
$$
|A+A|^2|AA|\log|A|\ge p|A|^2
$$
or
$$
|A+A|^2|AA|\log |A|\gg p|A|^5.
$$
In particular,
$$
\max\{|A+A|, |AA|\}\gg \min\Bigl\{|A|^{5/3}p^{-1/3}(\log
|A|)^{-1/3}, |A|^{2/3}p^{1/3} (\log|A|)^{-1/3}\Bigr\}.
$$

Thus, if we prove Lemma~\ref{lem:XYG}, then we are done. To this
end, let $I$ denote the number of solutions of the equation
$$
x+gy=x_1+gy_1,\quad g\in G, \,\, x, x_1\in X, \,\, y, y_1\in Y.
$$
We express $I$ via trigonometric sums:
\begin{eqnarray*}
& &I=\frac{1}{p}\sum_{n=0}^{p-1}\sum_{g\in G}\Bigl|\sum_{x\in
X}e^{2\pi i nx/p}\Bigr|^2\Bigl|\sum_{y\in Y}e^{2\pi i
ngy/p}\Bigr|^2.
\end{eqnarray*}
Picking up the term corresponding to $n=0$ and then extending the
summation over $g$ to the whole range $0\le g\le p-1,$ we obtain
\begin{eqnarray*} & & I
\le
\frac{|X|^2|Y|^2|G|}{p}+\frac{1}{p}\sum_{n=1}^{p-1}\sum_{g=0}^{p-1}\Bigl|\sum_{x\in
X}e^{2\pi i nx/p}\Bigr|^2\Bigl|\sum_{y\in Y}e^{2\pi i
ngy/p}\Bigr|^2\le\\
&&
\frac{|X|^2|Y|^2|G|}{p}+\frac{1}{p}\left(\sum_{n=0}^{p-1}\Bigl|\sum_{x\in
X}e^{2\pi i
nx/p}\Bigr|^2\right)\left(\sum_{g=0}^{p-1}\Bigl|\sum_{y\in
Y}e^{2\pi i gy/p}\Bigr|^2\right)= \\
&& \frac{|X|^2|Y|^2|G|}{p}+p|X||Y|.
\end{eqnarray*}
Therefore, there exists a fixed element $\xi\in G$ such that the
number $I_0$ of solutions of the equation
$$
x+\xi y=x_1+\xi y_1, \quad x,x_1\in X,\,\, y, y_1\in Y,
$$
satisfies
$$
I_0\le \frac{|X|^2|Y|^2}{p}+\frac{p|X||Y|}{|G|}.
$$
From the relationship between the cardinality of a set and the
number of solutions of the corresponding equation, we derive that
$$
|X+\xi*Y|\ge\frac{|X|^2|Y|^2}{I_0}\ge
\frac{|X|^2|Y|^2}{\frac{|X|^2|Y|^2}{p}+\frac{p|X||Y|}{|G|}}=\frac{p|X||Y||G|}{|X||Y||G|+p^2}.
$$
This proves Lemma~\ref{lem:XYG} and thus
Theorem~\ref{thm:sumprodest1} is also proved.

\bigskip

{\bf Acknowledgement.} The author is grateful to S.~V.~Konyagin and
I.~Z.~Ruzsa for very useful remarks, comments and references. This
work was supported by the Project PAPIIT IN 100307 from UNAM.

\end{document}